\documentclass[12pt]{amsart}

\def\squarebox#1{\hbox to #1{\hfill\vbox to #1{\vfill}}}

\newcommand{\Up}{\Upsilon}

\newcommand{\ep}{\epsilon}
\newcommand{\al}{\alpha}
\newcommand{\be}{\beta}

\newcommand{\cO}{{\mathcal O}}

\newcommand{\cV}{{\mathcal V}}
\newcommand{\cT}{{\mathcal T}}
\newcommand{\gd}{{\dot g}}

\newcommand{\R}{{\mathbb R}}
\newcommand{\Sb}{{\mathbb S}}
\newcommand{\N}{{\mathbb N}}

\newcommand{\Vol}{\operatorname{Vol}}
\newcommand{\Ric}{\operatorname{Ric}}

\theoremstyle{plain}
\newtheorem{theorem}{Theorem}[section]
\newtheorem{lemma}[theorem]{Lemma}

\theoremstyle{definition}

\theoremstyle{remark}

\numberwithin{equation}{section}

\title{The Ambient Obstruction Tensor and $Q$-Curvature}

\author{C. Robin Graham}
\address{Department of Mathematics, University of Washington,
Box 354350\\
Seattle, WA 98195-4350}
\email{robin@math.washington.edu}

\author{Kengo Hirachi}
\address{Graduate School of Mathematical Sciences, University of Tokyo\\
  3-8-1 Komaba, Megro, Tokyo 153-8914, Japan}
\email{hirachi@ms.u-tokyo.ac.jp}

\begin{document}
\maketitle

\thispagestyle{empty}

\renewcommand{\thefootnote}{}
\footnotetext{The work of the first author was partially supported by NSF 
grant DMS-0204480.  The work of the second author was supported by
Grant-in-Aid for Scientific Research, JSPS.  The authors are grateful to
the Erwin Schr\"odinger Institute for hospitality during the writing of
this paper.} 
\renewcommand{\thefootnote}{1}

\section{Introduction}\label{intro}
The Bach tensor is a basic object in four-dimensional conformal geometry.  
It is a conformally invariant trace-free symmetric 2-tensor involving 4
derivatives of 
the metric which is of particular interest because it vanishes for metrics
which are conformal to Einstein metrics, and because it arises as the first
variational derivative of the conformally invariant Lagrangian 
$\int |W|^2$, 
where $W$ denotes the Weyl tensor.  A generalization of the Bach tensor
to higher even dimensional manifolds was indicated in \cite{FG1}. 
This ``ambient obstruction tensor'', which, suitably normalized, we denote
by $\cO_{ij}$, is also  
a trace-free symmetric 2-tensor which is conformally invariant
and vanishes for conformally Einstein metrics.  It involves $n$
derivatives of the metric on a manifold of even dimension $n\geq 4$.
In this paper we give the details of the derivation and basic properties of   
the obstruction tensor and provide a characterization generalizing the
variational 
characterization of the Bach tensor in four dimensions.  We also give an 
invariant-theoretic classification of conformally invariant tensors up to
invariants which are quadratic and higher in curvature which illuminates
the fundamental nature of the obstruction tensor.

Our higher dimensional substitute for $\int |W|^2$ is the integral of
Branson's $Q$-curvature (\cite{B}).  This $Q$-curvature is a scalar
quantity defined on even-dimensional 
Riemannian (or pseudo-Riemannian) manifolds.  It is not a pointwise
conformal invariant like $|W|^2$, but it does have a simple transformation
law under conformal change which implies that its integral over a compact
manifold is a conformal
invariant.  In dimension 4, one has $6Q = -\Delta R +R^2 -3|\Ric|^2$,
where $R$ denotes the scalar curvature and $\Delta = \nabla^i\nabla_i$.   
Since the Pfaffian in 4 dimensions is a multiple of 
$R^2 -3|\Ric|^2 + \frac32|W|^2$, it follows that $\int Q $ is a linear
combination of the Euler characteristic and $\int |W|^2$, so 
the variational derivatives of $\int Q$ and $\int |W|^2$ are multiples of
one another.  It follows from a result announced by Alexakis
\cite{Al} that also in higher dimensions, $\int Q$ is a linear combination
of the Euler characteristic and the integral of a pointwise conformal
invariant.  However explicit formulae are not available.

Our variational characterization is:
\begin{theorem}\label{variation}
If $g^t$ is a 1-parameter family of metrics on a compact manifold $M$ of
even dimension $n\geq 4$, then
$$
\left( \int_M Q\,dv\right ){\dot {}}
= (-1)^{n/2}\,\frac{n-2}2 \int_M \cO_{ij}\gd^{ij} dv,
$$
where ${\dot {}} = \partial_t|_{t=0}$ and $\cO_{ij}$ and $dv$ on the right
hand side are with respect to $g^0$. 
\end{theorem}

In \cite{FG1}, $\cO_{ij}$ arose as the obstruction to the existence of a 
smooth formal power series solution for the ambient metric associated to
the given conformal 
structure, a Ricci-flat metric in 2 higher dimensions homogeneous with
respect to dilations.  As described in \cite{FG1}, the ambient metric is
equivalent 
to a Poincar\'e metric, a metric in 1 higher dimension with constant
negative Ricci curvature having the given conformal structure as conformal  
infinity, and the obstruction tensor may alternately be viewed as
obstructing smooth formal power series solutions for a Poincar\'e metric.  
It is the latter formulation that we use in this paper, for it is the
Poincar\'e metric that provides the link between $Q$-curvature and the
obstruction tensor.  Specifically, we use the result of \cite{GZ} 
that the integral of the
$Q$-curvature is equal to a multiple of the log term in the volume 
expansion of a Poincar\'e metric.  We then calculate the variation of the
log term coefficient by a simplified version of the method of Anderson
\cite{An} for expressing the variation of volume as a boundary integral.  
A different calculation of the first variation of the log coefficient in
the volume expansion is given in \cite{HSS}.  

The existence of the obstruction tensor gives rise to the questions of 
whether there are other conformally invariant tensors lurking in the
shadows, and whether there is some kind of odd-dimensional
analogue.  Of course, one may construct further invariants from known ones
by taking tensor products and contracting.  However, the following result
shows that up to quadratic and higher terms in curvature, the Weyl tensor
(or Cotton tensor in dimension 3) and the obstruction tensor are the
only irreducible conformally invariant tensors.
\begin{theorem}\label{invtensors}
A conformally invariant irreducible natural tensor of $n$-dimensional  
oriented Riemannian manifolds is equivalent modulo a conformally invariant 
natural tensor of degree at least $2$ in curvature with a multiple of one
of the following:  
\begin{itemize}
\item $n=3$:  the Cotton tensor $C_{ijk} = P_{ij},_k - P_{ik},_j$
\item $n=4$:  the self-dual or anti-self dual Weyl tensor $W^{\pm}_{ijkl}$
  or the Bach tensor $B_{ij}=\cO_{ij}$
\item $n\geq 5$ odd:  the Weyl tensor $W_{ijkl}$
\item $n\geq 6$ even: the Weyl tensor $W_{ijkl}$ or the obstruction tensor
  $\cO_{ij}$ 
\end{itemize}
\end{theorem}
\noindent
Here the trace-modified Ricci tensor $P_{ij}$ is defined by
\begin{equation}\label{P}
(n-2)P_{ij}=R_{ij}-\frac{R}{2(n-1)}g_{ij},
\end{equation}
and the terminology used in the statement of the Theorem is explained in \S 
4. 
Theorem~\ref{invtensors} is an easy consequence of the classification
of conformally invariant linear differential operators on the sphere due to
Boe-Collingwood (\cite{BC}).  

In \S 2 we show how to derive the obstruction tensor in terms of a
Poincar\'e metric and establish its basic properties.  In \S 3 we prove
Theorem~\ref{variation} and in \S 4 we prove Theorem~\ref{invtensors}.

We are grateful to Mike Eastwood and Andi \u{C}ap for helpful discussions
concerning the Boe-Collingwood classification of invariant operators.

\section{The Obstruction Tensor}\label{obst}
In this section we provide the details of the background
about the obstruction tensor.
We show how it arises 
as the obstruction to the existence of a smooth formal power series
solution for 
a Poincar\'e metric associated to the given conformal structure and derive
its properties from this characterization.    
We also show that this definition may be reformulated in terms  
of a formal solution to one higher order involving a log term.  

Let $M$ be a manifold of dimension $n\geq 3$ with smooth conformal
structure $[g]$ of 
signature $(p,q)$ and let $X$ be $n+1$-manifold with boundary $M$. 
All our considerations in this section are local near a point of $M$.   We
are interested 
in conformally compact metrics $g_+$ of signature $(p+1,q)$ on $X$ with 
conformal infinity $[g]$.
This means that if $x$ is a smooth defining
function for $M$, then $x^2 g_+$ is a smooth (to some order)
metric on $X$ with $x^2 g_+|_{TM} \in [g]$.  If $n$ is odd, then for any
conformal class $[g]$ there are metrics $g_+$ with $x^2 g_+$ a formal
smooth power 
series such that $\Ric g_+ = -ng_+$ to infinite order.  However if $n$ is
even, the obstruction tensor obstructs the existence of formal smooth
solutions at order $n-2$.  
(Throughout this paper, when we say that a tensor
is $O(x^s)$, we mean that all components of the tensor are $O(x^s)$ in a
smooth coordinate system on $X$.)\footnote{One could alternately consider
metrics of signature $(p,q+1)$ for which $\Ric g_+ = ng_+$.  This
formulation is equivalent to ours via the change $g_+\rightarrow -g_+$.} 
\begin{theorem}\label{obstensor}
If $n\geq 4$ is even, there exists a metric $g_+$ with $x^2g_+$ smooth such
that $g_+$ has $[g]$ as conformal infinity and 
$\Ric g_+ + ng_+ = O(x^{n-2})$.  $g_+$ is unique modulo $O(x^{n-2})$ up 
to a diffeomorphism of $X$ which restricts to the identity on $M$.  The
tensor $\operatorname{ tf } (x^{2-n}(\Ric g_+ + ng_+)|_{TM})$ 
on $M$ is independent of the choice of such $g_+$, where 
$\operatorname{tf}$
denotes the trace-free part with respect to $[g]$.  We define the
obstruction tensor  
\begin{equation}\label{defO}
\cO = c_n \operatorname{ tf } (x^{2-n}(\Ric g_+ + ng_+)|_{TM}),\qquad
c_n=\frac{2^{n-2}(n/2-1)!^2}{n-2}.
\end{equation}
Then $\cO_{ij}$ has the properties:
\begin{enumerate}
\item
$\cO$ is a natural tensor invariant of the metric $g=x^2g_+|_{TM}$;
i.e. in local coordinates the components of $\cO$ are given by
universal polynomials in the components of $g$, $g^{-1}$ and the curvature 
tensor of $g$ and its covariant derivatives.
The expression for $\cO_{ij}$  takes the form 
\begin{equation}\label{Oform}
\begin{split}
\cO_{ij} &= \Delta^{n/2-2}\left (P_{ij},_k{}^k-P_k{}^k,_{ij}\right )
+ lots\\
&=(3-n)^{-1}\Delta^{n/2-2}W_{kijl},{}^{kl} + lots,
\end{split}
\end{equation}
where  
$\Delta = \nabla^i\nabla_i$, $W$ denotes the Weyl tensor of $g$, 
and lots denotes quadratic and higher terms in curvature involving fewer 
derivatives.  
\item
One has
$$
\cO_i{}^i=0 \qquad\qquad\qquad \cO_{ij},{}^j = 0.
$$
\item
$\cO_{ij}$ is conformally invariant of weight $2-n$; i.e. if $0<\Omega \in 
  C^{\infty}(M)$ and $\hat{g}_{ij} =
\Omega^2 g_{ij}$, then $\hat{\cO}_{ij} = \Omega^{2-n}\cO_{ij}$.  
\item
If $g_{ij}$ is conformal to an Einstein metric, then $\cO_{ij}=0$.  
\end{enumerate}
\end{theorem}
\begin{proof}
There are discussions of the asymptotics of Poincar\'e metrics in 
the literature, but we provide a self-contained treatment.

We shall work with metrics in a normal form.  Lemma 5.2 and the subsequent
paragraph in \cite{GL} imply that if one is given a conformally compact
metric $g_+$ 
which is asymptotically Einstein in the sense that $\Ric g_+ + ng_+ =
O(x^{-1})$ and a representative metric $g \in [g]$, there is an
identification of a neighborhood of $M$ in $X$ with $M \times [0,\epsilon)$
such that $g_+$ takes the form 
\begin{equation}\label{specialform}
g_+ = x^{-2}(dx^2 + g_x)
\end{equation}
for a 1-parameter family $g_x$ of metrics on $M$ with $g_0 = g$.  

It is straightforward to calculate $E= \Ric g_+ + ng_+$ 
for $g_+$ of the form \eqref{specialform}.  
We use Greek indices to label objects on $X$,
Latin indices for $M$, and $0$ for $\partial_x$ so that in an
identification $X\cong M\times [0,\epsilon)$ as above, a Greek index
$\alpha$ corresponds to a pair $(i,0)$.  One obtains:
\begin{equation}\label{einstein1}
2xE_{ij}=-xg_{ij}'' + x g^{kl}g_{ik}'g_{jl}' - \frac{x}{2}
g^{kl}g_{kl}'g_{ij}'+(n-1)g_{ij}' 
+  g^{kl}g_{kl}'g_{ij} + 2x\Ric(g_x)_{ij}
\end{equation}
\begin{equation}\label{einstein2}
E_{i0} = \frac12 g^{kl}(\nabla_lg_{ik}'-\nabla_ig_{kl}')
\end{equation}
\begin{equation}\label{einstein3}
E_{00} = -\frac12g^{kl}g_{kl}''+ \frac14g^{kl}g^{pq}g_{kp}'g_{lq}' 
+ \frac12 x^{-1}g^{kl}g_{kl}',
\end{equation}
where $'$ denotes $\partial_x$, we have suppressed the subscript on $g_x$,   
and $\nabla$ and $\Ric$ denote the Levi-Citiva connection and Ricci
curvature of $g_x$ for fixed $x$.

One can determine the derivatives of $g_x$ inductively to solve the
equation $E_{ij}=O(x^{n-2})$, beginning with the prescription 
$g_0 = g$. Differentiating \eqref{einstein1} $s-1$ times 
and setting $x=0$ gives 
\begin{equation}\label{derivs}
\begin{split}
\partial_x^{s-1} (2xE_{ij})|_{x=0} = 
(n-s) & \partial_x^s g_{ij} + g^{kl}\partial_x^sg_{kl}g_{ij}\\
&+ \left ( \mbox{terms involving $\partial_x^kg_{ij}$ with $k< s$}
  \right ) 
\end{split}
\end{equation}
For $s \neq n,\, 2n$, the operator $\eta_{ij}\rightarrow (n-s)\eta_{ij} +
g^{kl}\eta_{kl}g_{ij}$ is invertible on symmetric 2-tensors at each point
of $M$.    It follows inductively that one uniquely
obtains a metric $g_+  \mod O(x^{n-2})$ 
of the form \eqref{specialform} by the requirement $E_{ij} = O(x^{n-2})$.
Moreover, the derivatives of $g_x$ at $x=0$ of order less than $n$ are all
natural tensor invariants of the initial representative metric $g$.  

The vanishing of the remaining components of $E$ to the correct order 
is deduced via the Bianchi identity. 
The Bianchi identity for Ricci curvature of $g_+$ states
$$
g_+^{\alpha\beta}\nabla^+_{\gamma}E_{\alpha\beta} =
2g_+^{\alpha\beta}\nabla^+_{\alpha}E_{\beta\gamma},
$$
where $\nabla^+$ denotes the Levi-Civita connection of $g_+$.  Taking 
separately $\gamma = 0$ and $\gamma = i$ and writing this
in terms of the connection $\nabla$ of $g_x$ gives the following two
equations:
\begin{equation}\label{bianchi1}
g^{jk}E_{jk}' = 2\nabla^jE_{j0} + (\partial_x + g^{jk}g_{jk}' 
- 2(n-1)x^{-1})E_{00}
\end{equation}
\begin{equation}\label{bianchi2}
\partial_iE_{00} + \nabla_iE_j{}^j - 2\nabla^jE_{ij}
= 2(\partial_x +\frac12g^{jk}g_{jk}' -(n-1) x^{-1})E_{i0}.
\end{equation}
We claim that $E_{00}=O(x^{n-2})$ and $E_{i0}=O(x^{n-1})$.  
These follow by induction on the statement that  
$E_{00}=O(x^{s-1})$ and $E_{i0}=O(x^{s})$ for
$0\leq s \leq n-1$.  The case $s=0$ is immediate from \eqref{einstein3} 
and \eqref{einstein2}.  Suppose the statement is true for some $s$, 
$s\leq n-2$.  Write $E_{00}=\lambda x^{s-1}$ and recall
$E_{ij}=O(x^{n-2})$.  In \eqref{bianchi1}, we 
have $g^{jk}E_{jk}' = O(x^{n-3}) = O(x^{s-1})$, 
$\nabla^jE_{j0} = O(x^s)$, and $g^{jk}g_{jk}'=O(1)$.  Calculating
\eqref{bianchi1} 
mod $O(x^{s-1})$ thus gives $(s-2n+1)\lambda x^{s-2} = O(x^{s-1})$, which 
implies $\lambda = O(x)$ so $E_{00}=O(x^s)$ as desired.  Now write
$E_{i0} = \mu_i x^s$ and calculate \eqref{bianchi2} mod $O(x^s)$.  
One obtains similarly $(s+1-n)\mu_i x^{s-1} = O(x^s)$.  Since $s\leq n-2$
it follows that $\mu_i = O(x)$ so $E_{i0} = O(x^{s+1})$, completing the
induction. 

This proves the first sentence of the statement of Theorem~\ref{obstensor}: 
existence of a formal solution to order $n-2$.    
The second sentence, uniqueness of $g_+$ up to diffeomorphism, 
follows from the fact that any metric $g_+$ can be put into the form
\eqref{specialform} by a diffeomorphism together with the uniqueness of the
determination of $g_+$ in the form \eqref{specialform} as above.
The definition \eqref{defO} of $\cO$ depends on a choice of defining
function $x$ to first order, equivalently on the choice of a conformal
representative,  but is otherwise diffeomorphism invariant.  So in order to 
establish the 
independence of $\cO$ on the freedom in $g_+$ at order $n-2$ and the
naturality of $\cO$, it suffices to consider $g_+$ of the form
\eqref{specialform}.   
These conclusions now follow from \eqref{derivs}:  taking $s=n$ shows that   
the trace-free part of $x^{2-n}E_{ij}|_{x=0}$ is given purely in terms of
the previously determined terms.  

The tensor $\cO_{ij}$ is trace-free by definition.  The fact that
$\cO_{ij},^j=0$ can be established by consideration of \eqref{bianchi1}
and \eqref{bianchi2} as follows.  Consider the approximately Einstein
metric $g_+ \mod O(x^{n-2})$ constructed inductively above which satisfies 
$E_{ij}=O(x^{n-2})$, $E_{i0}=O(x^{n-1})$, $E_{00}=O(x^{n-2})$.  
Although vanishing of the trace-free part of $E_{ij}$ at order $n-2$ 
is obstructed by $\cO_{ij}$, one sees from \eqref{derivs} that one can
solve for the trace to ensure that $g^{ij}E_{ij} = O(x^{n-1})$.  This is
sufficient to allow one to conclude exactly as above from \eqref{bianchi1}  
that $E_{00}=O(x^{n-1})$.  Now consider \eqref{bianchi2}.  One finds this
time that the right hand side is already $O(x^{n-1})$.  Substituting
$E_{ij} = \cO_{ij}x^{n-2} \mod O(x^{n-1})$ and calculating mod $O(x^{n-1})$
gives $\cO_{ij},^j=0$ as desired. 

The conformal invariance of $\cO_{ij}$ follows immediately from its
definition:  the rescaling $\hat{g}_{ij} =\Omega^2 g_{ij}$ corresponds to 
$\hat{x} = \Omega x +O(x^2)$, which by \eqref{defO} gives 
$\hat{\cO}_{ij} = \Omega^{2-n}\cO_{ij}$.  The fact that the
same tensor arises when calculated in the normal forms determined by
different conformal representatives is implicit in the invariance of the 
definition under diffeomorphisms.

The vanishing of $\cO_{ij}$ for conformally Einstein metrics follows from
the fact that for $g$ Einstein, one can write down an explicit solution  
for $g_+$.  It is well known that if $\Ric(g) = 4\lambda (n-1) g$, then the
metric 
$g_+ = x^{-2}(dx^2 + (1-\lambda x^2)^2 g)$ satisfies $\Ric(g_+) = -ng_+$. 
This is also easily checked directly using
\eqref{einstein1}--\eqref{einstein3}.  In particular there is no
obstruction to existence of a smooth formal
solution at order $n-2$, so it must be that $\cO_{ij}=0$.

To finish the proof of Theorem~\ref{obstensor}, it remains to derive the 
principal part of $\cO_{ij}$.  This can be done by keeping track of the
leading term in the inductive derivation above. As described above, the
derivatives $\partial_x^s (g_x)|_{x=0}$ for $1\leq s \leq n-1$ are
determined inductively by setting $E_{ij}=0$ and differentiating 
in \eqref{einstein1}, and
the obstruction $\cO_{ij}$ arises when trying to solve for 
$\partial_x^n(g_x)|_{x=0}$.  Parity considerations show that
these derivatives vanish for $s$ odd.  Differentiating \eqref{einstein1} once
gives $g_{ij}''|_{x=0} = -2P_{ij}$.
Differentiating further and using the first variation of Ricci curvature
$$
{\dot \Ric_{ij}} = \frac12({\dot g}_{ik},_j{}^k + {\dot g}_{jk},_i{}^k -
{\dot g}_{ij},_k{}^k - {\dot g}_k{}^k,_{ij})
$$
and the Bianchi identity $P_{ik},{}^k = P_k{}^k,_i$,
one determines inductively that 
$$
\partial_x^{2m}g_{ij}|_{x=0} = 2\frac{3\cdot5\cdot7\cdots(2m-1)}{(n-4)(n-6)
\cdots(n-2m)}
\left (\Delta^{m-2}P_k{}^k,{}_{ij}-\Delta^{m-1}P_{ij} \right ) 
+ \mbox{ lots }
$$
for $2\leq m<n/2$.  Using
$E_{ij} = c_n^{-1}x^{n-2}\cO_{ij} \mod O(x^{n-1})$ in \eqref{einstein1} and 
differentiating $n-1$ times then gives the first line of \eqref{Oform}.
The second follows from the fact that
$
W_{kijl},{}^{kl} = (3-n)(P_{ij},_k{}^k - P_{ik},_j{}^k).
$
\end{proof}

The proof of Theorem~\ref{obstensor} gives an algorithm for the calculation
of $\cO_{ij}$.  For $n=4,\,6$, carrying out the calculations gives the following
explicit formulae.  Define the Cotton and Bach tensors by: 
$$
C_{ijk} = P_{ij},_k - P_{ik},_j \qquad\qquad B_{ij} = P_{ij},_k{}^k 
- P_{ik},_j{}^k - P^{kl}W_{kijl}.
$$
Then when $n=4$ one has
$\cO_{ij} = B_{ij}$
and when $n=6$ one has
\[
\begin{split}
\cO_{ij} = B_{ij},_k{}^k - 2W_{kijl}&B^{kl} -4P_k{}^kB_{ij} +8P^{kl}C_{(ij)k},_l 
-4C^k{}_i{}^lC_{ljk} \\
&+2C_i{}^{kl}C_{jkl} +4P^k{}_k,_lC_{(ij)}{}^l
-4W_{kijl}P^k{}_mP^{ml}.
\end{split}
\]

If the obstruction tensor is nonzero, there are no formal smooth
solutions to $\Ric(g_+) = -n g_+$ beyond order $n-2$.  However, it is
always possible to find solutions to all orders by including log terms in
the expansion of $g_+$.  For our purposes it will suffice to consider
solutions to one higher order.  The obstruction tensor $\cO_{ij}$  
can then be characterized as the coefficient of the first log term.

\begin{theorem}\label{log}
In the setting of Theorem~\ref{obstensor}, there is a solution $g_+$ 
to $\Ric(g_+) + ng_+ = O(x^{n-1}\log x)$ of the form 
$g_+ = x^{-2}(dx^2 + g_x)$, where $g_x = h_x + r_x x^n\log x$ and 
$h_x$ and $r_x$ are smooth in $x$.  The coefficient $r_x$ is uniquely
determined at $x=0$ and is given by $nc_nr_0= 2 \cO$.  
\end{theorem}
\begin{proof}
Fix a metric $g_+^0=x^{-2}(dx^2 + g_x^0)$ with $g_x^0$ smooth 
which solves $E^{0}= O(x^{n-2})$ as
in Theorem~\ref{obstensor}. Set $g_x = g_x^0 + rx^n\log x +s x^n$.  
Substituting into \eqref{einstein1} gives
\[
\begin{split}
2xE_{ij}=2xE_{ij}^0 + g^{kl}r_{kl}g_{ij}(&x^{n-1}\log x +x^{n-1})\\
&-nr_{ij}x^{n-1} +ng^{kl}s_{kl}g_{ij}x^{n-1} \mod O(x^n\log x).
\end{split}
\]
It is required that this expression vanish mod $O(x^n\log x)$.
The requirement that there be no $x^{n-1}\log x$ term forces
$g^{kl}r_{kl}|_{x=0}=0$.  Since 
$c_n \operatorname{ tf }(x^{2-n} E_{ij}^0)|_{x=0} =\cO_{ij}$, we must 
have $nc_nr_{ij}|_{x=0} = 2 \cO_{ij}$.  The trace of $s_{ij}$ can be
chosen to guarantee $g^{ij}(x^{2-n}E_{ij}|_{x=0})=0$; the trace-free part  
of $s_{ij}$ remains arbitrary.  With these choices, if we set
$h_x = g_x^0 +sx^n \mod O(x^{n+1})$ and $r_x = r \mod O(1)$, we obtain  
$g_+$ in the form required in the statement of the theorem satisfying 
$E_{ij}=O(x^{n-1}\log x)$.  

In the proof of Theorem~\ref{obstensor} it was shown that $g_+^0$ satisfies  
$E_{i0}^0=O(x^{n-1})$, $E_{00}^0=O(x^{n-2})$.
It is evident from this and \eqref{einstein2}, \eqref{einstein3} that  
$g_+$ satisfies $E_{i0}=O(x^{n-1}\log x)$, $E_{00}=O(x^{n-2}\log x)$. 
Arguing as 
in the proof of Theorem~\ref{obstensor}, one finds that \eqref{bianchi1}  
implies that in fact one has $E_{00}=O(x^{n-1}\log x)$, 
completing the proof.
\end{proof}

\section{Proof of Theorem~\ref{variation}}
For any metric $g$ on $M$, 
we consider a metric $g_+ = x^{-2}(dx^2 + g_x)$ on $M\times
(0,\epsilon)$ given by Theorem~\ref{log} which satisfies
$\Ric(g_+) + ng_+ = O(x^{n-1}\log x)$.  Recall that $g_x$
was determined only up to addition of a trace-free tensor in the
coefficient of $x^n$ and up to addition of terms of order greater
than $n$.  For definiteness, we specify $g_x$ to be given
by the finite expansion
\begin{equation}\label{expansion}
g_x = g + g^{(2)}x^2 + \mbox{ (even powers) } + \frac2{nc_n} 
\cO x^n \log x +g^{(n)}x^n
\end{equation}
where the $g^{(2m)}$ for $m<n/2$ are those coefficients derived in 
Theorem~\ref{obstensor},  and we take $g^{(n)}$ to be the multiple of $g$
determined in the proof of Theorem~\ref{log}.  Then the metric
$g_+$ is completely determined by $g$, and of course satisfies 
$\Ric (g_+) + ng_+ = O(x^{n-1}\log x)$.  

The proof of Theorem~\ref{variation} depends on the volume expansion of an
asymptotically Einstein metric; see \cite{G}.  
The volume form of $g_+$ is 
$$
dv_{g_+} = x^{-n-1}\left( \frac{\det g_x}{\det g} \right )^{1/2}dv_gdx.
$$
From \eqref{expansion} and the fact that $g^{ij}\cO_{ij}=0$, it follows
that 
\begin{equation}\label{detexp}
\left( \frac{\det g_x}{\det g} \right )^{1/2}
= 1 +v^{(2)}x^2+ \mbox{ (even powers) } + v^{(n)}x^n +\,\, \cdots,
\end{equation}
where the $v^{(2j)}$ are locally determined invariant scalars given 
in terms of $g$ and its curvature, and $\cdots$ denotes terms vanishing to
higher order.  Integrating, it follows that for fixed $\ep_0$ we have the
asymptotic expansion as $\ep\rightarrow 0$
$$
\Vol_{g_+}(\{\ep<x<\ep_0\}) = c_0\ep^{-n}+c_2 \ep^{-n+2}
+\mbox{ (even powers) } +c_{n-2}\ep^{-2} +L\log{\frac1\ep} +O(1),
$$
where $c_{2j}=(n-2j)^{-1}\int_Mv^{(2j)}dv_g$ and $L=\int_Mv^{(n)}dv_g$.
The log term coefficient $L$ is invariant under conformal rescalings of
$g$; see \cite{G} for a proof.

The $Q$-curvature of a metric $g$ was originally defined by Branson
\cite{B} in terms of the zero-th order term of the conformally invariant
$n$-th power of the Laplacian $P_n$ of \cite{GJMS} by dimensional
continuation.  It is  
a scalar quantity with a particularly simple 
transformation law under conformal rescalings:  if ${\hat g} = e^{2\Up} g$,  
then $e^{n\Up} {\hat Q} = Q + P_n\Up$.  Although
$Q$ is not pointwise conformally invariant, it follows from the facts that 
$P_n$ is self-adjoint and annihilates constants that the 
integral $\int Q\,dv$ over a compact manifold is a conformal
invariant.  
Characterizations of the $Q$-curvature in terms of Poincar\'e metrics 
were given in \cite{GZ} and \cite{FG2}, and in terms of the ambient metric
in \cite{FH}.  We refer to \cite{B} and these
references for background about $Q$-curvature.  The main fact we will need
here is the result of \cite{GZ} (or see \cite{FG2} for a simpler proof)
that 
\begin{equation}\label{LQ}
\int Q\,dv = k_n L,\qquad  k_n = (-1)^{n/2}n(n-2)c_n.
\end{equation}

According to \eqref{LQ}, in order to calculate the variation of 
$\int Q\,dv$ it suffices to compute $\dot L$.  For this, we
use a simplification of the method of Anderson \cite{An} to rewrite the
variation of volume as a 
boundary integral for variations through Einstein metrics with fixed scalar
curvature.  In our case we need to estimate the errors resulting from the
fact that our metrics are only asymptotically Einstein.
\begin{lemma}\label{volvar}
Let $g^t$ be a 1-parameter family of metrics on a compact manifold $M$ and 
let $g_+^t = x^{-2}(dx^2 + g_x^t)$ be the corresponding asymptotically
Einstein metrics on $M\times (0,\ep_0)$, where for each $t$, $g_x^t$ is
constructed from $g^t$ as in \eqref{expansion}.  Set $X_{\ep} =
\{\ep<x<\ep_0\}$.  Then as $\ep\rightarrow 0$
we have
\begin{equation}\label{volform}
%\begin{split}
\Vol_{g_+^t} (X_{\ep}){\dot {}}\\ 
= \frac{\ep^{1-n}}{2n} 
\int_{x=\ep} \left (-\frac12 g^{ij}g^{kl}g_{jl}'{\dot g_{ik}} 
+ x^{-1}g^{ij}{\dot g_{ij}} - (g^{ij}{\dot g_{ij}})'\right ) dv_{g_{\ep}}
+O(1).
%\end{split}
\end{equation}
On the right hand side, $'$ denotes $\partial_x$ and ${\dot {}}$ denotes 
$\partial_t|_{t=0}$ as usual, and all $g=g_x^t$ are evaluated at $x=\ep$,
$t=0$ (after differentiation).  
\end{lemma}
\begin{proof}
Since the metrics $g_+^t$ are asymptotically Einstein, we
have uniformly in $t$ (and supressing the $t$-dependence of $g_+^t$):  
$$
\Ric_{g_+} = -ng_+ + O(x^{n-1}\log x),
$$
$$
R_{g_+} = -n(n+1) +O(x^{n+1}\log x),
$$
$$
{\dot {R_{g_+}}} = O(x^{n+1}\log x).
$$
Therefore $\int_{X_{\ep}}{\dot {R_{g_+}}} dv_{g_+} = O(1)$ as 
$\ep \rightarrow 0$.

On the other hand, the usual formula for the first variation of scalar
curvature gives
\[
\begin{split}
{\dot {R_{g_+}}}&= {\dot {g_+}}_{\al\be},^{\al\be} -{\dot
  {g_+}}_{\al}{}^{\al},_{\be}{}^{\be} - 
\Ric_{g_+}^{\al\be}{\dot {g_+}}_{\al\be} \\
&=
{\dot {g_+}}_{\al\be},^{\al\be} -{\dot
  {g_+}}_{\al}{}^{\al},_{\be}{}^{\be}  
+n g_+^{\al\be}{\dot {g_+}}_{\al\be} + O(x^{n+1}\log x),
\end{split}
\]
where the covariant derivatives are with respect to 
the Levi-Civita connection of $g_+$ and indices are raised and lowered
using $g_+$.  Integrating gives
$$
\int_{X_{\ep}}({\dot {g_+}}_{\al\be},^{\al\be}
-{\dot {g_+}}_{\al}{}^{\al},_{\be}{}^{\be})dv_{g_+}  
+2n\int_{X_{\ep}}{\dot {dv_{g_+}}} = O(1),
$$
so
\[
\begin{split}
-2n\Vol_{g_+} (X_{\ep}){\dot {}}&=
\int_{X_{\ep}}({\dot {g_+}}_{\al\be},^{\al\be}
-{\dot {g_+}}_{\al}{}^{\al},_{\be}{}^{\be})dv_{g_+}  +O(1)\\
&=\int_{\partial X_{\ep}}({\dot {g_+}}_{\al\be},^{\al}
-{\dot {g_+}}_{\al}{}^{\al},_{\be})\nu_+^{\be}d\sigma_+ +O(1),
\end{split}
\]
where $\nu_+^{\be}$
denotes the outward unit normal and $d\sigma_+$ the  
induced volume density.  The integral over $x=\ep_0$ is $O(1)$, and on  
$x=\ep$ we have  $\nu_+^{\be} = -\ep \delta^{\be}_0$,
$d\sigma_+ = \ep^{-n}dv_{g_{\ep}}$.  Also,
${\dot {g_+}}_{\al\be}$ vanishes if either $\al=0$ or $\be=0$, 
and ${\dot {g_+}}_{ij} = x^{-2}\gd_{ij}$.  An easy computation
relating the connections of $g_+$ and $g_x$ shows that 
$$
{\dot {g_+}}_{\al 0},^{\al}-{\dot {g_+}}_{\al}{}^{\al},_0
=-\frac12 g^{ij}g^{kl}g_{jl}'{\dot g_{ik}} 
+ x^{-1}g^{ij}{\dot g_{ij}} - (g^{ij}{\dot g_{ij}})',
$$
which gives \eqref{volform}.
\end{proof}

Now ${\dot L}$ occurs as the coefficient of $\log \frac1{\ep}$ in the
asymptotic expansion of the left hand side of \eqref{volform}.  So we need
to evaluate the  
$\ep^{n-1}\log \frac1{\ep}$ coefficient in the expansion of the integral on
the right hand
side.  From $dv_{g_x} = (\det g_x/\det g)^{1/2}\,dv_g$ and 
\eqref{detexp}, it follows that the expansion of the volume form 
has no $x^n\log x$ term, so does not contribute to this coefficient. 
Differentiating the volume form in $t$, one concludes that also the
expansion of 
$g^{ij}\gd_{ij}$ has no $x^n\log x$ term.  Therefore the second and
third terms in the integrand also do not contribute to the 
$\ep^{n-1}\log \frac1{\ep}$ coefficient.  The
only contribution from the first term in the integrand comes from the log
term in $g_{jl}'$.  This gives $2nc_n{\dot L} =
\int_M\cO_{ij}\gd^{ij}\,dv_g$, which combined with \eqref{LQ} gives
Theorem~\ref{variation}.  

\section{Proof of Theorem~\ref{invtensors}}
We consider natural tensor invariants of oriented $n$-dimensional
Riemannian manifolds $(M,g)$ with values in a subbundle $\cV\subset
\otimes^k T^*M$ induced by 
a representation of $SO(n)$ on an invariant subspace 
$V\subset \otimes^k (\R^n)^*$; see \cite{E} for a discussion of natural
tensors.  The components of a natural tensor are expressible as
linear combinations of partial contractions of the metric,
the volume form, the Riemannian curvature tensor, and 
its covariant derivatives.  
We say that a natural tensor is irreducible if it is
nonzero and if $V$ is irreducible as an $SO(n)$-module.  We say that two
natural tensors with values in subbundles $\cV_1$, $\cV_2$ are equivalent
if they correspond under an isomorphism $\cV_1\cong \cV_2$ induced by an
$SO(n)$-module isomorphism $V_1\cong V_2$ of the underlying subspaces.  

The Ricci identity for commuting covariant derivatives does not preserve
homogeneity degree, so the degree of a natural tensor as a polynomial in
curvature and its derivatives is not well-defined.  However, the space of
natural tensors is filtered by degree and it does make sense to say that a
natural tensor is of degree at least $d$ in curvature for $d\in \N$.
A natural tensor $\cT(g)$ is said to be conformally invariant of weight $w$
if  $\cT(\Omega^2 g) = \Omega^w \cT(g)$ for $0<\Omega \in C^{\infty}(M)$.
The naturality of $\cT$ implies that if $\phi$ is any local diffeomorphism,
then $\cT(\phi^*g)=\phi^*\cT(g)$.

The idea of the proof of Theorem~\ref{invtensors} is simple:  linearizing a
natural tensor at the usual  
metric on the sphere $\Sb^n$ gives a linear differential operator on
infinitesimal  
metrics, and the conformal invariance implies that this differential
operator satisfies an invariance property under conformal motions.  A known
theorem classifies such invariant differential operators and this
classification in the linear case implies the classification up to
quadratic and higher terms for natural tensors. 

In more detail, let $\cT(g)$ be an irreducible natural tensor which is 
conformally invariant of weight $w$.  Denote by $g_0$ the usual metric on
$\Sb^n$.  Define a linear differential operator $T$ from the bundle 
$S^2_0T^*\Sb^n$ of trace-free symmetric 2-tensors on $\Sb^n$ to the bundle
$\cV$ on $\Sb^n$ by 
$$
T(h) = \frac{d}{dt} \cT(g_0 +th)|_{t=0}
$$
for $h$ a section of $S^2_0T^*\Sb^n$.  We claim that if 
$\phi:\Sb^n\rightarrow \Sb^n$ is a conformal motion of $\Sb^n$ satisfying 
$\phi^*g_0 = \Omega^2 g_0$ for $0<\Omega \in C^{\infty}(\Sb^n)$, then
\begin{equation}\label{invariance}
\Omega^w T(\phi^*h)=\phi^*( T(\Omega^2\circ\phi^{-1}\,h)).
\end{equation}
In fact, 
\[
\begin{split}
\Omega^w \cT(g_0 +t\phi^*h)& = \Omega^w \cT(\Omega^{-2}\phi^*g_0 +t\phi^*h)
= \Omega^w \cT(\Omega^{-2}\phi^*(g_0 +t\Omega^2\circ\phi^{-1}\,h))\\
&= \cT(\phi^*(g_0 +t\Omega^2\circ\phi^{-1}\,h)) 
= \phi^*\cT(g_0 +t\Omega^2\circ\phi^{-1}\,h),
\end{split}
\]
from which the claim follows by differentiation.  

Let $G=O_e(n+1,1)$
denote the identity component of the conformal group and 
$P\subset G$ the isotropy group of a point on $\Sb^n$, so that $\Sb^n=G/P$.
Then  
\eqref{invariance} states exactly that $T$ is a $G$-equivariant map between   
sections of the homogeneous bundles $S^2_0T^*\Sb^n(2)$ and $\cV(w)$ on
$G/P$, where the number in parentheses indicates the conformal weight
of the homogeneous bundle.
Such invariant differential operators between any irreducible homogeneous 
bundles on $\Sb^n$ have been completely classified; see (1.4) of \cite{BC}.
The classification in \cite{BC} is formulated in terms of the homomorphisms
of the generalized Verma modules dual to the homomorphisms of the modules
of jets of sections of the homogeneous bundles induced by the 
differential operators.  See \cite{BE}, \cite{ES} and references cited
there for elaboration and interpretation of this classification in the 
context of conformal geometry.  

For our purposes it is sufficient to know all the invariant
operators with domain $S^2_0T^*\Sb^n(2)$ and range any irreducible bundle.  
The bundle $S^2_0T^*\Sb^n(2)$ has regular integral infinitesimal character
so fits into a generalized Bernstein-Gelfand-Gelfand complex of
invariant operators (the so called deformation complex -- see \cite{GG} for
a direct construction of this complex on a general conformally flat
manifold).  In odd 
dimensions, up to scale and equivalence, there is precisely one invariant
operator acting on $S^2_0T^*\Sb^n(2)$: the linearized Cotton 
tensor in dimension 3 and the linearized Weyl tensor in higher dimensions.  
In even dimensions there are further operators.  For $n\geq 6$  
there is one more operator acting on $S^2_0T^*\Sb^n(2)$, which must be 
the linearized obstruction tensor since this is a nonzero operator acting
between the appropriate bundles.  The case $n=4$ is exceptional 
because $S^2_0T^*\Sb^n(2)$ occurs at the edge of the middle diamond of the
Hasse diagram, and there are three invariant operators acting on
$S^2_0T^*\Sb^n(2)$, which can be identified with the linearized self-dual
and anti-self-dual Weyl tensors and the linearized Bach tensor.  

Theorem~\ref{invtensors} is an immediate consequence.  The linearization 
of a conformally invariant irreducible 
natural tensor is equivalent to a multiple of one of the
invariant operators given by the Boe-Collingwood classification.  By
inspection, each such 
operator is the linearization of one of the conformally invariant natural
tensors listed in the statement of the Theorem.  
But if two conformally invariant natural tensors have the same
linearization on the sphere, their 
difference must be of degree at least 2 in curvature.

\end{document}